\documentclass[3p]{elsarticle}

\usepackage{amsmath}
\usepackage{amsthm}
\usepackage{amssymb}
\usepackage{bm}
\usepackage{accents}

\newtheorem{thm}{Theorem}

\newdefinition{rmk}{Remark}
\newproof{pf}{Proof}

\numberwithin{equation}{section}

\graphicspath{{figs/}}

\journal{arXiv}

\begin{document}

\begin{frontmatter}

\title{Computational decomposition and composition technique for approximate solution of nonstationary problems\tnoteref{label1}}
\tnotetext[label1]{The work was supported by the Russian Science Foundation (grant No. 23-41-00037).}

\author{P.N. Vabishchevich\corref{cor1}\fnref{lab1,lab2}}
\ead{vab@cs.msu.ru}
\cortext[cor1]{Correspondibg author.}

\address[lab1]{Lomonosov Moscow State University, 1, building 52, Leninskie Gory,  119991 Moscow, Russia}

\address[lab2]{North-Eastern Federal University, 58, Belinskogo st, Yakutsk, 677000, Russia}

\begin{abstract}
Stable computational algorithms for the approximate solution of the Cauchy problem for nonstationary problems are based on implicit time approximations. Computational costs for boundary value problems for systems of coupled multidimensional equations can be reduced by additive decomposition of the problem operator(s) and composition of the approximate solution using particular explicit-implicit time approximations.
Such a technique is currently applied in conditions where the decomposition step is uncomplicated.
A general approach is proposed to construct decomposition-composition algorithms for evolution equations in finite-dimensional Hilbert spaces.
It is based on two main variants of the decomposition of the unit operator in the corresponding spaces at the decomposition stage and the application of additive operator-difference schemes at the composition stage.
The general results are illustrated on the boundary value problem for a second-order parabolic equation by constructing standard splitting schemes on spatial variables and region-additive schemes (domain decomposition schemes).
\end{abstract}

\end{frontmatter}

\section{Introduction}\label{sec:1}

Mathematical modeling is based on numerical solutions of boundary value problems for systems of multidimensional partial equations.
Computational algorithms are usually based on finite-element and finite-volume approximations over the space \cite{KnabnerAngermann2003,QuarteroniValli}.
Implicit time approximations are often used when numerically solving nonstationary boundary value problems for partial derivative equations. Such approximations provide unconditional stability of the solution concerning the initial data and the right-hand side \cite{Ascher2008,LeVeque2007}.
Explicit schemes, which are simpler to find an approximate solution on a new level in time, have tight constraints on the grid spacing in time \cite{Samarskii1989,SamarskiiMatusVabischevich2002}.
We want to construct schemes that are as stable as implicit schemes and as easy to implement computationally as explicit schemes.

When inhomogeneous time approximations are used, the problem operator is split into two operator summands with the allocation of a computationally acceptable summand, which is taken from the upper level in time and the other summand --- from the lower level.
Such IMEX methods are widely used in computational practice\cite{Ascher1995,HundsdorferVerwer2003}.
For example, in \cite{Vabishchevich2020}, explicit- implicit two- and three-level operator-difference schemes are constructed for the first-order evolution equation for both the standard splitting of the primary operator of the problem and the splitting of the operator at the time derivative of the solution.

Splitting schemes {\cite{Marchuk1990,VabishchevichAdditive} are based on the known additive representation of the problem operator.
In this case, the transition to a new level in time is carried out by solving evolutionary problems for individual operator summands.
In many nonstationary problems, computationally acceptable subproblems make sense when constructed on the principle of solution decomposition when more straightforward problems are formulated for separate components of the solution.
Such new solution splitting schemes are constructed in \cite{efendiev2021splitting,vabishchevich2021solution} for the approximate solution of the Cauchy problem in finite-dimensional Hilbert space for evolution equations of first and second order.

We can consider splitting schemes for approximate solutions of nonstationary problems as a computational technology for the decomposition (analysis) of the situation and composition (synthesis) of the solution.
At the stage of decomposition, an additive representation of the problem's operator(s) into more straightforward operators is performed, and at the stage of composition --- an approximate solution of the problem is constructed from the solutions of the issues for individual operator summands on the basis of special time approximations.
Various classes of two- and three-level additive operator-difference schemes (splitting schemes) {\cite{VabishchevichAdditive} under a given additive splitting of the problem operator have been constructed.
Explicit-implicit approximations are most easily constructed for two-component splitting.
We use component-wise splitting schemes (sum approximation schemes), regularized additive schemes, and vector splitting schemes for general multi-component splitting.

In computational practice, most attention is paid to splitting schemes on spatial variables or physical processes when decomposing the problem operator(s) is not difficult.
We develop a more general approach to construct an additive decomposition of the problem operators in the approximate solution of the Cauchy problem for evolution equations in Hilbert finite-dimensional spaces.
The key idea is related to using an additive representation of the unit operator in appropriate spaces.
In the composition stage, additive operator-difference schemes are used.
From these positions, we can consider, in particular, the previously proposed region-additive schemes (domain decomposition schemes) \cite{vabishchevich2023subdomain,vabishchevich2023difference} for nonstationary problems, which is based on the partitioning of the unit for the domain.

The paper is organized as follows.
Section \ref{sec:2} describes the problem of constructing splitting schemes for approximate solutions of nonstationary issues based on the decomposition and composition techniques.
We consider the Cauchy problem in a finite-dimensional Hilbert space for a first-order evolution equation with a self-adjoint operator.
The model boundary value problem for a parabolic equation in a rectangle under finite-difference approximation is considered in Section \ref{sec:3}.
The decomposition is provided by splitting by spatial variables (directional variable schemes) and based on the computational grid decomposition (region-additive schemes).
The decomposition operators are constructed in Section \ref{sec:4} section.
The first approach (decompositions in one space) is based on the additive representation of the unit operator in the corresponding space as a sum of self-adjoint operators. When investigating the stability of the composition schemes, the operators of the corresponding evolution equations are symmetrized.
In Section \ref{sec:5}, we formulate a second approach to the decomposition stage when working on a set of spaces.
The decomposition operators are represented in factorized form, which allows us to construct an additive splitting not only for  problem operator but also for the solution itself.
Section \ref{sec:6} is devoted to briefly discussing the decomposition-composition technique for other problems.
In particular, the problems for second-order evolution equations and systems of first-order equations are considered.
The results are summarized in Section \ref{sec:7}.

\section{Problem statement}\label{sec:2}

Additive operator schemes (splitting schemes) are considered in the general computational technology of decomposition-composition.
In the example of the Cauchy problem for the first-order evolution equation, the problem of constructing an additive splitting of the problem operator is formulated.
Splitting schemes for the two-dimensional parabolic equation using a uniform rectangular grid over the space are given as a benchmark.

\subsection{Cauchy problem}

Let us consider the Cauchy problem for the first-order evolution equation:
\begin{equation}\label{2.1}
\frac{d u}{d t} + A u = f(t),
\quad 0 < t \leq T,
\end{equation}
\begin{equation}\label{2.2}
u(0)= u^0 .
\end{equation}
We are looking for a solution $u(t)$ of the equation \eqref{2.1} for $0 < t \leq T$ from a finite-dimensional Hilbert space $H$ under a given initial condition \eqref{2.2}.
The constant (independent of $t$) operator $A$ in (\ref{2.1}) is self-adjoint and positive:
\begin{equation}\label{2.3}
A = A^* > 0 .
\end{equation}

In $H$ the scalar product for $u, v \in H$ is $(u,v)$ and the norm is $\|u\| = (u,u)^{1/2}$.
For a self-adjoint and positive operator $D$, we define a Hilbert space
$H_D$ with scalar product and norm $(u,v)_D = (D u,v), \ \|u\|_D = (u,v)_D^{1/2}$.
For the solution of the problem (\ref{2.1})--(\ref{2.3}), the simplest a priori estimates are as follows
\begin{equation}\label{2.4}
\|u(t)\|_D^2 \leq \|u^0\|_D^2 + \frac{1}{2} \int_{0}^{t} \|f(s)\|^2_{D A^{-1}} \, d s,
\quad 0 < t \leq T,
\end{equation}
in $H_D$ when, for example, $D = A, I, A^{-1}$.
Such estimates of the stability of the solution on the initial data and the right-hand side should be inherited by using some time approximations.

We introduce a uniform, for simplicity, time grid with step $\tau$ and let $y^n = y(t^n), \ t^n = n \tau$, $n = 0,1, \ldots, N, \ N\tau = T$.
To numerically solve the problem (\ref{2.1}), (\ref{2.2}), we will use a two-level scheme with weights ($\sigma = \mathrm{const}$):
\begin{equation}\label{2.5}
\frac{y^{n+1} - y^{n}}{\tau } + A (\sigma y^{n+1} + (1-\sigma) y^{n}) = f^{n+\sigma} ,
\quad n = 0,1, \ldots, N-1 ,
\end{equation}
\begin{equation}\label{2.6}
y^0 = u^0 ,
\end{equation}
when using the notation
\[
t^{n+\sigma} = \sigma t^{n+1} + (1-\sigma) t^n,
\quad f^{n+\sigma} = \sigma f^{n+1} + (1-\sigma) f^n .
\]
The difference scheme (\ref{2.5}), (\ref{2.6}) approximates the problem (\ref{2.1}), (\ref{2.2}) with second order $\tau$ at $\sigma = 0.5$ (symmetric scheme) and with first order --- at $\sigma \neq 0.5$.

When formulating stability conditions for two- and three-level schemes, we refer to general stability results for operator-difference schemes \cite{Samarskii1989,SamarskiiMatusVabischevich2002}.

\begin{thm}\label{t-1}
The scheme with weight (\ref{2.5}), (\ref{2.6}) is unconditionally stable at $\sigma \geq 1/2$ in $H_D$, $D = A, I, A^{-1}$.
For the solution, the estimate
\begin{equation}\label{2.7}
\|y^{n+1}\|_D^2 \leq \| u^0 \|_D^2 + \frac{1}{2} \sum_{k=0}^{n} \tau \| f^{k+\sigma}\|_{DA^{-1}}^2,
\quad \ n = 0,1, \ldots, N-1 ,
\end{equation}
is holds.
\end{thm}

\begin{pf}
We write (\ref{2.5}) as
\begin{equation}\label{2.8}
(I + \tau G) \frac{y^{n+1} - y^{n}}{\tau } + A \frac{y^{n+1} + y^{n}}{2} = f^{n+\sigma} ,
\quad \ n = 0,1, \ldots, N-1 ,
\end{equation}
with the operator
\[
G = \Big (\sigma - \frac{1}{2} \Big ) A .
\]
At $\sigma \geq 1/2$ we have $G = G^* \geq 0$.
By multiplying (\ref{2.8}) scalarly in $H$ by $2 D A^{-1}(y^{n+1} - y^{n})$, we get
\[
\begin{split}
2 \tau \Big \| \frac{y^{n+1} - y^{n}}{\tau } \Big \|_{DA^{-1}}^2 & + 2 \big (DA^{-1}G(y^{n+1} - y^{n}), y^{n+1} - y^{n} \big )
+ \|y^{n+1}\|_{D}^2 \\
& = \|y^{n}\|_{D}^2 + 2 \tau \Big (DA^{-1} f^{n+\sigma} , \frac{y^{n+1} - y^{n}}{\tau } \Big ) .
\end{split}
\]
Considering
\[
\begin{split}
\big (DA^{-1}G(y^{n+1} - y^{n}, y^{n+1} - y^{n} \big ) & \geq 0, \\
\Big (DA^{-1} f^{n+\sigma} , \frac{y^{n+1} - y^{n}}{\tau } \Big ) & \leq \Big \| \frac{y^{n+1} - y^{n}}{\tau } \Big \|_{DA^{-1}}^2
+ \frac{1}{4} \| f^{n+\sigma}\|_{DA^{-1}}^2,
\end{split}
\]
we obtain the inequality
\[
\|y^{n+1}\|_D^2 \leq \|y^{n}\|_D^2 + \frac{1}{2} \tau \| f^{n+\sigma}\|_{DA^{-1}}^2 .
\]
From this follows the provable estimate (\ref{2.7}).
\end{pf}

To find the solution on the new time level, we solve the problem
\[
(I + \sigma \tau A) y^{n+1} = \psi^{n},
\quad \psi^{n} = \big(I - (1-\sigma) \tau A \big) y^{n+1} + \tau f^{n+\sigma} ,
\quad \ n = 0,1, \ldots, N-1 .
\]
Several variants of iterative methods can be considered for this purpose.
The main feature of unsteady problems is that we have an excellent initial approximation $y^{n}$ for the unknown solution $y^{n+1}$.
This properties often limits the number of iterations.

\subsection{Decomposition-composition technique}

The general methodological technique of studying systems based on analysis and synthesis is used to reduce the computational effort in finding an approximate solution at a new time level.
Simpler private subproblems are selected (analysis), the study of which allows for obtaining solutions to the problem's general problem (synthesis).
Such computational technology of analysis (decomposition) and synthesis (composition) in solving unsteady problems is realized by using splitting schemes.
Let us consider the critical elements of this approach in the example of the Cauchy problem (\ref{2.1}), (\ref{2.2}).

The difficulties in solving the problem are related to the complexity of the operator $A$.
We use the additive representation of the operator $A$ as a sum of simpler operators at the decomposition stage.
Let us assume that in our problem (\ref{2.1})--(\ref{2.3}) takes place
\begin{equation}\label{2.9}
A = \sum_{\alpha =1}^{p} A_\alpha ,
\quad A_\alpha = A_\alpha^* \geq 0,
\quad \alpha = 1,2, \ldots, p .
\end{equation}
According to the decomposition (\ref{2.9}), the solution of the auxiliary problems for the equations 
\begin{equation}\label{2.10}
\frac{d u_\alpha}{d t} + A_\alpha u_\alpha = f_\alpha(t)
\end{equation}
with separate operator summands $A_\alpha^* \geq 0, \ \alpha = 1,2, \ldots, p$ is a simpler problem than the solution of (\ref{2.1}), (\ref{2.2}) with operator $A$.

At the composition stage, the approximate solution of the initial problem $u(t)$ is constructed (synthesized) from the solutions of auxiliary problems with operators $A_\alpha^* \geq 0, \ \alpha = 1,2, \ldots, p$, which are determined from equations (\ref{2.10}).
For this purpose, we construct additive operator-difference schemes \cite{Samarskii1989,VabishchevichAdditive} when decomposing (\ref{2.9}).

The current state of research and practical use of the decomposition-composition technique is characterized by a relatively deep elaboration of the composition stage.
Different splitting schemes for the additive representation (\ref{2.9}) of the problem operator have been proposed, and unconditionally stable operator-difference schemes have been obtained.
Less attention has been paid to the problem of decomposition. This is mainly because additive splitting itself is often natural and unproblematic in the considered problems.
This paper aims to formulate general approaches for more straightforward problems in decomposing a nonstationary problem.

\section{Model problem}{\label{sec:3}}

We illustrate the possibilities of constructing splitting schemes on the example of a problem for a second-order parabolic equation in a rectangle.
Standard difference approximations on a uniform rectangular grid are used for discretization over space.
The decomposition is provided by splitting the operator over the spatial variables.
The second case is associated with the decomposition of the domain into subdomains.
Various additive operator-difference schemes are used in the composition.

\subsection{Differential Problem}

We consider a model boundary value problem for a second-order parabolic equation in a rectangle
\[
\Omega = \{ {\bm x} \ | \ {\bm x} = (x_1, x_2),
\ 0 < x_{\alpha} < l_{\alpha}, \ \alpha =1,2 \}.
\]
The unknown function $w({\bm x},t)$ satisfies the equation
\begin{equation}\label{3.1}
\frac{\partial w}{\partial t}
- \sum_{\alpha =1}^{2}
\frac{\partial }{\partial x_\alpha}
\left ( k({\bm x}) \frac{\partial w}{\partial x_\alpha} \right ) = \varphi ({\bm x},t),
\quad {\bm x}\in \Omega,
\quad 0 < t \leq T,
\end{equation}
in which $k({\bm x}) \geq \kappa > 0, \ {\bm x}\in \Omega$.
Let us supplement equation (\ref{3.1}) with homogeneous Dirichlet boundary conditions
\begin{equation}\label{3.2}
w({\bm x},t) = 0,
\quad {\bm x}\in \partial \Omega,
\quad 0 < t \leq T.
\end{equation}
The initial condition is also set
\begin{equation}\label{3.3}
w({\bm x},0) = w^0({\bm x}),
\quad {\bm x}\in \Omega.
\end{equation}

The nonstationary diffusion problem (\ref{3.1}), (\ref{3.3}) is considered on the set of functions $w({\bm x},t)$ satisfying the boundary conditions (\ref{3.2}).
Then instead of (\ref{3.1}), (\ref{3.2}) we can use the differential-operator equation
\begin{equation}\label{3.4}
\frac {dw}{dt} + {\cal A} u = \varphi(t),
\quad 0 < t \leq T.
\end{equation}
We consider the Cauchy problem for the evolution equation (\ref{3.4}):
\begin{equation}\label{3.5}
w(0) = w^0 .
\end{equation}
For the elliptic operator, we have
\[
{\cal A} w =
- \sum_{\alpha =1}^{2}
\frac{\partial }{\partial x_\alpha}
\left ( k({\bm x}) \frac{\partial w}{\partial x_\alpha} \right ) .
\]

On the set of sufficiently smooth functions (\ref{3.2}), define a Hilbert space ${\cal H} = {\cal L}_2 (\Omega)$ with scalar product and norm
\[
(w, v) = \int_{\Omega} w({\bm x}) v({\bm x}) d {\bm x},
\quad \| w\| = (w,w)^{1/2} .
\]
In ${\cal H}$, the operator ${\cal A}$ is self-adjoint and positively defined:
\begin{equation}\label{3.6}
{\cal A} = {\cal A}^* \ge \kappa \delta {\cal I},
\quad \delta = \delta(\Omega) > 0,
\end{equation}
where ${\cal I}$ is the unit operator in ${\cal H}$.

\begin{figure}[ht].
\center
\includegraphics[width=0.49\linewidth]{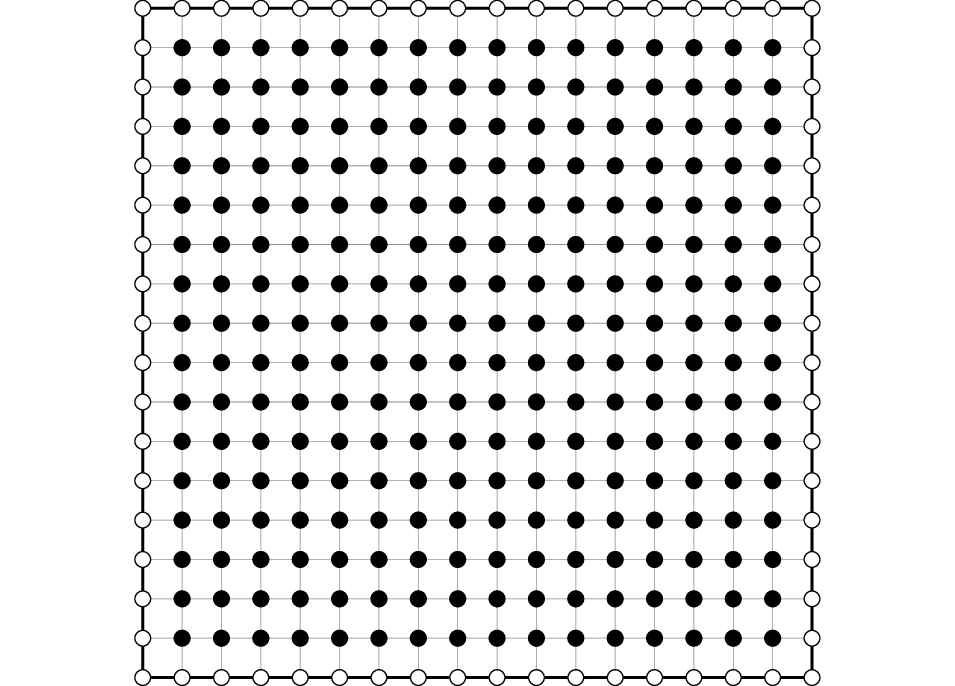}
\includegraphics[width=0.49\linewidth]{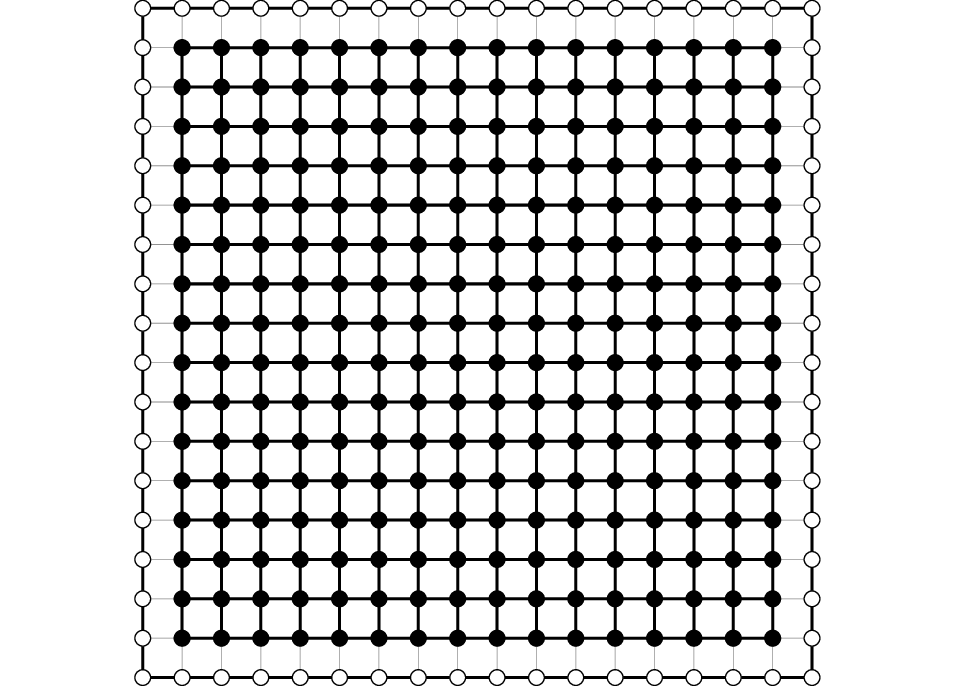}
\caption{Grid: $\bullet$ --- internal nodes, $\circ$ --- boundary nodes}.
\label{f-1}
\end{figure}

\subsection{Approximation over space}

The approximate solution is given at the nodes of a uniform rectangular grid in $\Omega$:
\[
\bar{\omega} = \{ {\bm x} \ | \ {\bm x} = (x_1, x_2),
\quad x_\alpha = i_\alpha h_\alpha,
\quad i_\alpha = 0,1,...,N_\alpha,
\quad N_\alpha h_\alpha = l_\alpha\}.
\]
and let $\omega$ be the set of internal nodes ($\bar{\omega} = \omega \cup \partial \omega$) (see Figure \ref{f-1}).
For the grid functions $u({\bm x}) = 0, \ {\bm x} \in \partial \omega$, we define a Hilbert space $H = L_2({\omega})$ with scalar product and norm
\[
(u,v) = \sum_{{\bm x} \in \omega}
u({\bm x}) v({\bm x}) h_1 h_2,
\quad \|u\| = (u,u)^{1/2} .
\]

Considering the coefficient $k({\bm x})$ in the domain of $\Omega$ smooth enough, we take the grid elliptic operator as
\begin{equation}\label{3.7}
\begin{split}
A u & =
- \frac{1}{h_1^2} k(x_1+0.5h_1,x_2)
(u(x_1+h_1,x_2) - u(x_1,x_2)) \\
& + \frac{1}{h_1^2} k(x_1-0.5h_1,x_2)
(u(x_1,x_2) - u(x_1-h_1,x_2)) \\
& - \frac{1}{h_2^2} k(x_1,x_2+0.5h_2)
(u(x_1,x_2+h_2) - u(x_1,x_2)) \\
& + \frac{1}{h_2^2} k(x_1,x_2-0.5h_2)
(u(x_1,x_2) - u(x_1,x_2-h_2)) .
\end{split}
\end{equation}
In $H$, the operator $A$ is self-adjoint and positively defined \cite{Samarskii1989}:
\begin{equation}\label{3.8}
A = A^* \geq \kappa (\delta_1+\delta_2) I,
\quad \delta_{\alpha} =
\frac{4}{h^2_{\alpha}} \sin^2 \frac{\pi h_{\alpha}}{2 l_{\alpha}} ,
\quad \alpha = 1,2.
\end{equation}

After space approximation from (\ref{3.4}), (\ref{3.5}) we arrive at the Cauchy problem (\ref{2.1}), (\ref{2.2}) when $u^0({\bm x}) = w^0({\bm x}), \ \bm x \in \omega$.
The involved nodal values in equation (\ref{2.1}) when using the approximation (\ref{3.7}), (\ref{3.8}) are labeled with lines in Fig.\ref{f-1} on the right.

\subsection{Decomposition}

The classical variant of splitting schemes is related to splitting the operator by spatial variables.
In the case of the problem (\ref{3.4})--(\ref{3.8}) we have a two-component ($p=2$) additive representation (\ref{2.9}) of the operator $A$ at
\[
\begin{split}
A_1 u & =
- \frac{1}{h_1^2} k(x_1+0.5h_1,x_2)
(u(x_1+h_1,x_2) - u(x_1,x_2)) \\
& + \frac{1}{h_1^2} k(x_1-0.5h_1,x_2)
(u(x_1,x_2) - u(x_1-h_1,x_2)) ,
\end{split}
\]
\[
\begin{split}
A_2 u & =
- \frac{1}{h_2^2} k(x_1,x_2+0.5h_2)
(u(x_1,x_2+h_2) - u(x_1,x_2)) \\
& + \frac{1}{h_2^2} k(x_1,x_2-0.5h_2)
(u(x_1,x_2) - u(x_1,x_2-h_2)) .
\end{split}
\]
The involved grid nodes for the operators $A_1$ and $A_2$ are shown in Fig.\ref{f-2}.

\begin{figure}[ht]
\center
\includegraphics[width=0.49\linewidth]{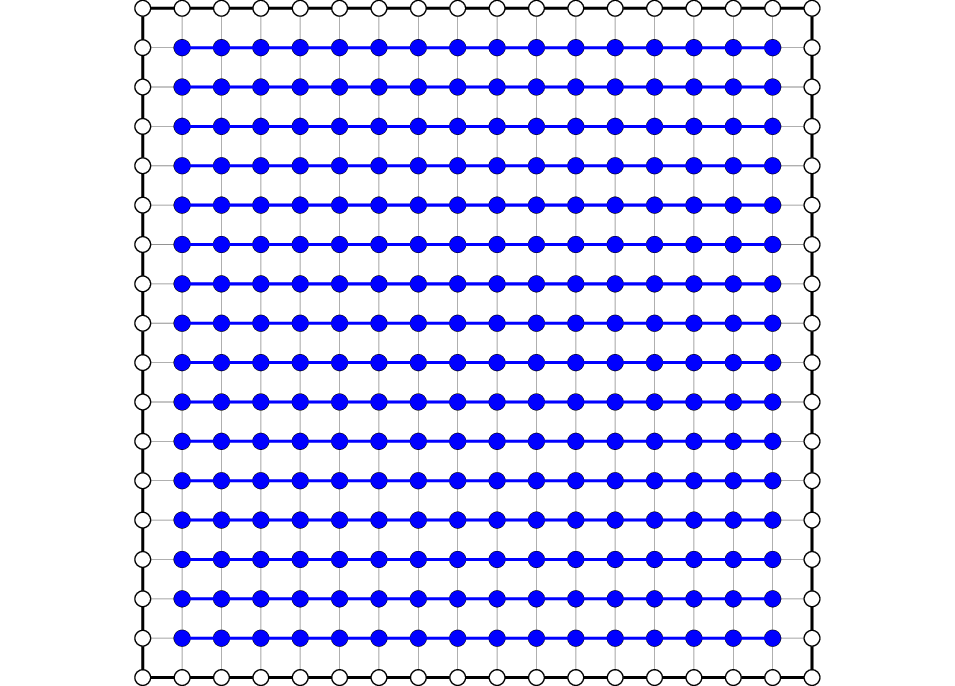}
\includegraphics[width=0.49\linewidth]{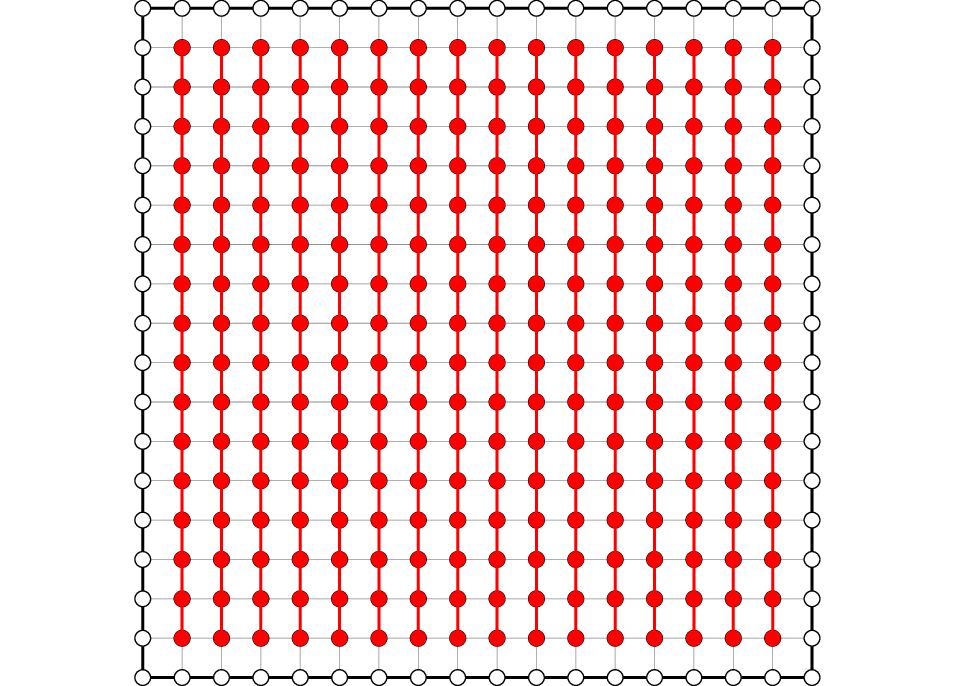}
\caption{Calculated nodes for $A_1$ (left) and for $A_2$ (right).}
\label{f-2}
\end{figure}

We also note a class of noniterative domain decomposition methods for approximate solution of multidimensional initial boundary value problems \cite{SamarskiiMatusVabischevich2002}.
In the approximate solution of the problem (\ref{3.1})--(\ref{3.3}) we will use the domain decomposition
\begin{equation}\label{3.9}
\overline{\Omega} = \bigcup_{\alpha =1}^{p} \overline{\Omega}_{\alpha},
\quad \overline{\Omega} = 
\Omega_{\alpha} \cup \partial \Omega_{\alpha},
\quad \alpha = 1, 2, ..., p ,
\end{equation}
with overlapping subdomains
($\Omega_{\alpha \beta} \equiv \Omega_{\alpha } \cap \Omega_{\beta} \neq \varnothing$)
or without overlapping subdomains ($\Omega_{\alpha \beta} = \varnothing$).
Based on (\ref{3.9}), we construct the corresponding additive representation of the problem operator
\begin{equation}\label{3.10}
\mathcal{A} = \sum_{\alpha =1}^{p} \mathcal{A}_{\alpha} .
\end{equation}
The operator $\mathcal{A}_{\alpha}$ is associated with the solution of some problem in the subdomain $\Omega_{\alpha}, \ \alpha = 1,2, ..., p$.

We will construct the operators $A_{\alpha}, \ \alpha = 1,2, ..., p$ using the unit partitioning for the computational domain.
In the decomposition (\ref{3.9}), we associate a function $\eta_{\alpha}({\bm x})$ with a separate subdomain $\Omega_{\alpha}$
($\alpha = 1,2,...,p$) such that
\begin{equation}\label{3.11}
\eta_{\alpha}({\bm x}) = \left \{
\begin{array}{cc}
> 0, & {\bm x} \in \Omega_{\alpha},\\
0, & {\bm x} \notin \Omega_{\alpha}, \\
\end{array}
\right .
\quad \alpha = 1,2,...,p ,
\end{equation}
provided that
\begin{equation}\label{3.12}
\sum_{\alpha =1}^{p} \eta_{\alpha}({\bm x}) = 1,
\quad {\bm x} \in \Omega .
\end{equation}
For the problem at hand (\ref{3.1})--(\ref{3.3}), we can use \cite{SamarskiiMatusVabischevich2002,vabishchevich1989difference} the following constructions:
\begin{equation}\label{3.13}
\mathcal{A}_{\alpha} = \eta_{\alpha} \, \mathcal{A},
\end{equation}
\begin{equation}\label{3.14}
\mathcal{A}_{\alpha} = \mathcal{A} \, \eta_{\alpha} ,
\quad \alpha = 1,2,...,p .
\end{equation}

We do the same to obtain an additive representation of the discrete operator $A$.
We associate the domain decomposition (\ref{3.9}) with distinct subsets of grid nodes $\omega_{\alpha}, \ \alpha =1,2, ..., p$:
\[
\omega = \mathop{\cup}^{p}_{\alpha=1} \omega_\alpha,
\quad \omega_{\alpha} = \{ {\bm x} \ | \ {\bm x} \in \omega, \ {\bm x} \in \Omega_\alpha \},
\quad \alpha =1,2, ..., p .
\]
We compare the partitioning of the unit domain (\ref{3.11}), (\ref{3.12}) to the partitioning for the set of internal nodes of $\omega$:
\begin{equation}\label{3.15}
\sum_{\alpha =1}^{p} \chi_{\alpha} = 1,
\quad \chi_{\alpha} (\bm x) \geq 0,
\quad \bm x \in \omega ,
\quad \alpha = 1,2,...,p .
\end{equation}
In the simplest case, $\chi_{\alpha} (\bm x) = \eta_{\alpha} (\bm x), \ \ \bm x \in \omega$.
Fig.\ref{f-3} illustrates the separation of grid parts for the case of non-intersecting two subdomains.

\begin{figure}[ht]
\center
\includegraphics[width=0.49\linewidth]{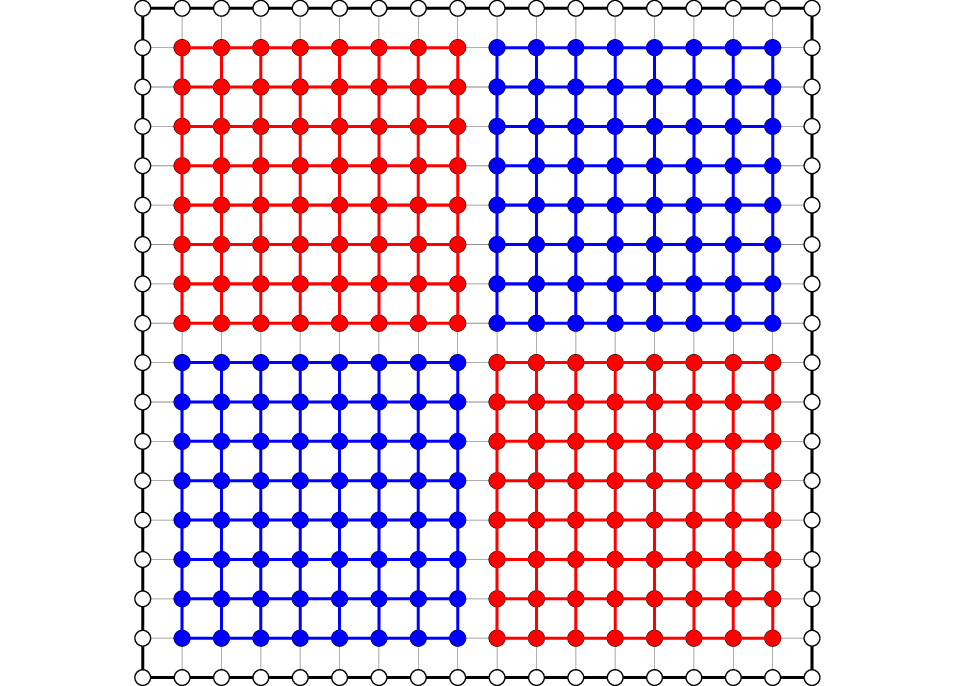}
\includegraphics[width=0.49\linewidth]{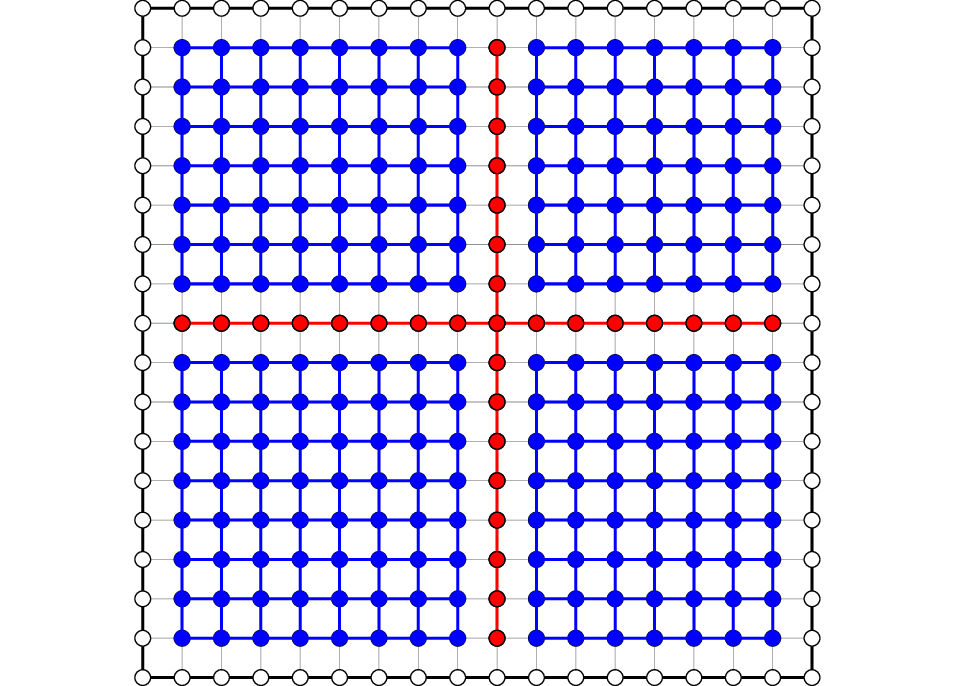}
\caption{Grid when the region is partitioned into rectangles (left) and when skeletonized (right): $\omega_1$ --- blue nodes, $\omega_2$ --- red nodes.}
\label{f-3}
\end{figure}

Similar to (\ref{3.13}), (\ref{3.14}), we define the decomposition operators as
\begin{equation}\label{3.16}
A_{\alpha} = \chi_{\alpha} A,
\end{equation}
\begin{equation}\label{3.17}
A_{\alpha} = A \chi_{\alpha} ,
\quad \alpha = 1,2,...,p .
\end{equation}
By (\ref{3.15}), when (\ref{3.16}), (\ref{3.17}) is split, there is an additive representation for the problem operator
\begin{equation}\label{3.18}
A = \sum_{\alpha =1}^{p}A_{\alpha} .
\end{equation}
The condition of non-negativity and self-adjointness of the operators $A_{\alpha}, \ \alpha = 1,2,...,p$ (see (\ref{2.9})) is provided by the symmetrization procedure \cite{SamarskiiMatusVabischevich2002}.
In our model problem (\ref{3.1})--(\ref{3.3}), neighboring nodes from other subsets of the grid are involved to compute the values of $A_{\alpha} u$ at nodes $\bm x \in \omega_\alpha , \ \alpha = 1,2,...,p$.

\subsection{Composition} 

In the composition stage, the approximate solution is constructed from solutions of auxiliary problems with separate operator summands (\ref{2.10}) in an additive representation (\ref{3.18}).
We consider the Cauchy problem for Eq.
\begin{equation}\label{3.19}
\frac{d u}{d t} + \sum_{\alpha =1}^{p}A_{\alpha} y = f(t)
\end{equation}
under the chosen additive representation of the right-hand side:
\[
f(t) = \sum_{\alpha =1}^{p}f_{\alpha}(t) .
\]
Various additive operator-difference schemes \cite{VabishchevichAdditive} are used to numerically solve the problem (\ref{2.2}), (\ref{3.19}).
Let us note the main classes of splitting schemes and formulate conditions for their unconditional stability.
In addition to the decomposition (\ref{2.9}), we can consider the more general case with non-self-conjugate operator summands:
\begin{equation}\label{3.20}
A = \sum_{\alpha =1}^{p} A_\alpha ,
\quad A_\alpha \geq 0,
\quad \alpha = 1,2, \ldots, p .
\end{equation}

Let us separately emphasize the case of two-component decomposition ($p=2$).
For the problems (\ref{2.2}), (\ref{3.19}), (\ref{3.20}) we can use factorized schemes \cite{Samarskii1989}
\begin{equation}\label{3.21}
(I + \sigma \tau A_1) (I + \sigma \tau A_2) \frac{y^{n+1} - y^{n}}{\tau } + (A_1 + A_2) y^{n} = f^{n+\sigma} ,
\quad n = 0,1, \ldots, N-1 .
\end{equation}
When the weight parameter $\sigma = 0.5$, we have the operator analog of the classical Peaceman--Rachford scheme \cite{PeacemanRachford1955}, and when $\sigma = 1$ --- the Douglas--Rachford scheme \cite{DouglasRachford1956}.
The factorized schemes (\ref{2.6}), (\ref{3.20}), (\ref{3.21}) are unconditionally stable in $H_D$ with $D = (I + \sigma \tau A_2^*) (I + \sigma \tau A_2)$ for $\sigma \geq 0.5$. For the approximate solution,  we have the estimate
\[
\big \| (I + \sigma \tau A_2) y^{n+1} \big \| \leq \big \| (I + \sigma \tau A_2) u^0 \big \| +
\sum_{k=0}^{n} \tau \big \| f^{n+\sigma} \big \| .
\]
When $\sigma = 0.5$, the factorized scheme (\ref{2.6}), (\ref{3.21}) has a second order accuracy of $\tau$, and when $\sigma \neq 0.5$ --- a first order accuracy.

For general multicomponent decomposition ($p > 1$ in (\ref{3.18})), we can consider the schemes of component-wise splitting \cite{Yanenko1959,Yanenko1967}.
The transition to a new level in time (solution composition) in the problem (\ref{2.2}), (\ref{3.19}) is provided by sequential solution of auxiliary problems
\begin{equation}\label{3.22}
\begin{split}
\frac{y^{n+\alpha/p} - y^{n+(\alpha-1)/p}} {\tau}
& + A_\alpha (\sigma y^{n+\alpha/p}
+ (1-\sigma)y^{n+(\alpha-1)/p}) \\
& = f^{n+\sigma}_\alpha,
\quad \alpha = 1,2,\ldots,p, \quad n = 0,1,\ldots, N-1.
\end{split}
\end{equation}
The stability of the scheme (\ref{2.6}), (\ref{3.20}), (\ref{3.22}) holds in $H$ under the usual weight constraints $\sigma \geq 0.5, \ \alpha = 1,2,\ldots,p$.
The study of convergence of the approximate solution at integer steps is based on the concept of sum approximation \cite{Samarskii1962}.
The component-wise splitting scheme (\ref{2.6}), (\ref{3.22}) has first-order $\tau$ accuracy.
When we choose $\sigma = 0.5, \ \alpha = 1,2,\ldots,p$ and organize the calculations according to the Fryazinov-Strang regulation \cite{Fryazinov1968,Strang1968}
\[
A_1 \rightarrow A_2 \rightarrow \cdots \rightarrow A_p \rightarrow A_p \rightarrow \cdots \rightarrow A_1 ,
\]
we can count on the second order of accuracy.

Let us highlight a variant of construction of additive difference schemes of component-wise splitting, in which independent solutions of simple problems and parallel organization of computations are allowed.
In additive-averaged schemes \cite{VabishchevichAdditive,GordezianiMeladze1974}, the transition to a new time level is carried out as follows:
\begin{equation}\label{3.23}
\begin{split}
\frac{y^{n+1}_\alpha - y^{n}} {p \tau} +
& A_\alpha (\sigma y^{n+1}_\alpha
+ (1-\sigma)y^{n}) = f^{n+\sigma}_\alpha, \\\
& \alpha = 1,2,...,p, \quad n = 0,1,\ldots, N-1 ,\\\\\
& y^{n+1} = \frac 1p \sum_{\alpha=1}^{p} y^{n+1}_\alpha .
\end{split}
\end{equation}
The composition of the solution at the new time level is provided by averaging the solutions of the auxiliary problems.
The stability conditions for the additive-averaging schemes (\ref{2.6}), (\ref{3.20}), (\ref{3.23}) are the same as for the conventional component-wise splitting schemes (\ref{2.6}), (\ref{3.20}), (\ref{3.22}): stability in $H$ takes place at $\sigma \geq 0.5$.

A more convenient construction of time approximations in multicomponent decomposition is realized by applying regularized additive schemes \cite{VabishchevichAdditive,SamarskiiVabischevich1998}.
It is based on multiplicative regularization of the explicit scheme, where the approximate solution is determined from Eq.
\begin{equation}\label{3.24}
\frac{y^{n+1} - y^{n}} {\tau} + \sum_{\alpha=1}^{p} (I + \sigma \tau A_\alpha )^{-1} A_\alpha y^{n}
= f^{n+\sigma} .
\end{equation}
The stability constraints on the time step in the explicit scheme are related to the norm of the problem operator.
The multiplicative regularization used aims to reduce this norm appropriately.
In the regularized additive scheme, no intermediate problems or auxiliary functions are introduced, and the original equation is approximated directly.
Sufficient conditions for unconditional stability of the scheme (\ref{2.6}), (\ref{3.20}), (\ref{3.24}) in $H$ are $\sigma \geq 0.5p$.

When constructing vector additive schemes \cite{Abrashin1990,VabishchevichVector} instead of the scalar equation (\ref{3.19}), we solve the Cauchy problem for a system of identical equations:
\begin{equation}\label{3.25}
\frac {d u_\alpha}{dt} +
\sum_{\beta=1}^{p}A_{\beta} u_{\beta} = f(t),
\quad t > 0,
\end{equation}
\begin{equation}\label{3.26}
u_\alpha(0) = u^0 ,
\quad \alpha=1,2,\ldots,p .
\end{equation}
In this case, $u_\alpha(t) = u(t), \ \alpha=1,2,\ldots,p$ and therefore any component of the vector $\bm u =\{u_1, u_2, \ldots, u_p\}$ can be taken as the solution to the original problem (\ref{2.2}), (\ref{3.19}).
We construct various variants of two- and three-level vector composition schemes \cite{VabishchevichAdditive} for the approximate solution of the problem (\ref{3.25}), (\ref{3.26}).
An example is the scheme
\[
(I + \sigma\tau A_\alpha)
\frac{y_\alpha^{n+1} - y_\alpha^{n}}{\tau} +
\sum_{\beta=1}^{p} A_\beta y_\beta^{n} =
\varphi^n,
\quad \alpha=1,2,\ldots,p ,
\]
which is stable at $\sigma \geq 0.5 p$.

\section{Decomposition in one space}\label{sec:4}

A sequence of more straightforward problems is constructed based on the decomposition of the problem operator.
We apply a general approach based on an additive representation of the unit operator.
We highlight various general constructions of such decomposition and note the possibilities of constructing composition schemes.

\subsection{Additive decomposition of the unit operator}

We construct time approximations for the approximate solution of the Cauchy problem (\ref{2.1})--(\ref{2.3}) based on the decomposition of the operator $A: H \rightarrow H$ in the form (\ref{3.18}). 
The operator summands $A_\alpha , \ \alpha =1,2, \ldots, p$ are obtained based on the general approach, without being directly connected to a particular form of the operator $A$.

The unit (identity) operator $I$ in $H$ can be represented as
\begin{equation}\label{4.1}
I = \sum_{\alpha = 1}^{p} R_\alpha ,
\quad R_\alpha = R_\alpha^* \geq 0,
\quad \alpha=1,2,\ldots,p .
\end{equation}
Depending on the applications, the operators $R_\alpha: H \rightarrow H, \ \alpha =1,2, \ldots, p$ can be associated with different variants of restriction operators.
In particular, we can mention the classical projection operators.
When applying a unit domain decomposition for the computational grid (\ref{3.15}), we have
\[
R_\alpha u = \chi_\alpha (\bm x) u,
\quad \bm x \in \omega ,
\quad \alpha=1,2,\ldots,p .
\]

Starting from the decomposition of the unit operator (\ref{4.1}) we can proceed to the decomposition (\ref{3.18}) of the operator $A$ in two main ways:
\begin{equation}\label{4.2}
A \longrightarrow \sum_{\alpha = 1}^{p} R_\alpha A,
\quad A_\alpha = R_\alpha A,
\end{equation}
\begin{equation}\label{4.3}
A \longrightarrow \sum_{\alpha = 1}^{p} A R_\alpha ,
\quad A_\alpha = A R_\alpha ,
\quad \alpha=1,2,\ldots,p .
\end{equation}
We noted similar constructions (see (\ref{3.16}), (\ref{3.17})) when considering domain decomposition schemes for the numerical solution of the model parabolic problem.

\subsection{The first variant of decomposition}

Let us consider the use of decompositions (\ref{4.2}) and (\ref{4.3}) on the example of the Cauchy problem (\ref{2.1})--(\ref{2.3}).
If we use the variant (\ref{4.2}), we obtain the equation
\begin{equation}\label{4.4}
\frac{d u}{d t} + \sum_{\alpha = 1}^{p} R_\alpha A u = f(t),
\quad 0 < t \leq T .
\end{equation}
Direct use of the results of the theory of additive operator-difference schemes for the problem (\ref{2.2}), (\ref{2.3}), (\ref{4.4}) is not possible because the operators $R_\alpha A, \ \alpha =1,2, \ldots, p$ are not non-negative.
The transformation of the problem by symmetrization of \cite{SamarskiiMatusVabischevich2002,SamVabGulin} is salvageable.

Let's multiply the equation by $A^{1/2}$ and get the equation
\begin{equation}\label{4.5}
\frac{d \widetilde{u}} {d t} + \sum_{\alpha = 1}^{p} \widetilde{A}_\alpha \widetilde{u} = \widetilde{f}(t),
\end{equation}
in which $\widetilde{u} = A^{1/2} u, \ \widetilde{f} = A^{1/2} f$ and
\[
\widetilde{A}_\alpha = A^{1/2} R_\alpha A^{1/2},
\quad \widetilde{A}_\alpha = \widetilde{A}_\alpha^* \geq 0,
\quad \alpha=1,2,\ldots,p .
\]
The initial condition (\ref{2.2}) takes the form of
\begin{equation}\label{4.6}
\widetilde{u} (0) = \widetilde{u}^0,
\quad \widetilde{u}^0 = A^{1/2} u^0 .
\end{equation}
Multiplying equation (\ref{4.5}) scalarly in $H$ by $\widetilde{u}$, we get
\begin{equation}\label{4.7}
\frac{1}{2} \frac{d }{d t} \| \widetilde{u} (t) \|^2 +
\sum_{\alpha = 1}^{p} \Big ( \widetilde{A}_\alpha \widetilde{u}, \widetilde{u} \Big ) = \big (\widetilde{f}(t),\widetilde{u} \big ) .
\end{equation}
Given (\ref{4.1}) and the notation introduced, we have
\[
\| \widetilde{u} (t) \|^2 = \|u(t)\|_A^2,
\quad \sum_{\alpha = 1}^{p} \Big ( \widetilde{A}_\alpha \widetilde{u}, \widetilde{u} \Big ) = \|A u\|^2,
\quad \big (\widetilde{f}(t),\widetilde{u} \big ) = (Au, f) .
\]
Inequality
\[
\frac{d }{d t} \|u(t)\|_A^2 \leq \frac{1}{2} \|f(t) \||^2
\]
follows from from (\ref{4.7}).
We obtain an a priori estimate (\ref{2.4}) for $D = A$.

For the problem (\ref{4.5}), (\ref{4.6}), we can use different time approximations.
We will not load the text of the paper with technical details that are not principal.
When discussing the problems of composition schemes, we will consider the homogeneous equation (\ref{2.1}) (right-hand side $f(t) = 0$).
For example, when applying the regularized additive scheme (see (\ref{3.24})) to approximate the solution of the problem (\ref{2.2}), (\ref{4.2}) we use the scheme
\begin{equation}\label{4.8}
\frac{y^{n+1} - y^n}{\tau} + \sum_{\alpha = 1}^{p} (I+\sigma \tau R_\alpha A)^{-1} R_\alpha A y^n = 0 .
\end{equation}

\begin{thm}\label{t-2}
The scheme (\ref{2.6}), (\ref{4.8}) is unconditionally stable at $\sigma \geq p/2$ in $H_A$.
For the approximate solution there is an estimation of stability on initial data
\begin{equation}\label{4.9}
\|y^{n+1}\|_A \leq \| u^0 \|_A,
\quad \ n = 0,1, \ldots, N-1 .
\end{equation}
\end{thm}

\begin{pf}
We write the equation (\ref{4.8}) as a system of equations
\begin{equation}\label{4.10}
\frac{y_\alpha^{n+1} - y^n}{p \tau} + (I+\sigma \tau R_\alpha A)^{-1} R_\alpha A y^n = 0 ,
\quad \alpha=1,2,\ldots,p ,
\end{equation}
when the solution is additively represented at the new level:
\begin{equation}\label{4.11}
y^{n+1} = \frac{1}{p} \sum_{\alpha = 1}^{p} y_\alpha^{n+1} .
\end{equation}
The organization of the computation of (\ref{4.10}), (\ref{4.11}) can be directly related to the use of additive-averaged splitting schemes (see (\ref{3.23})).
From (\ref{4.10}) with $\sigma \geq p/2$ we obtain
\[
\Big (I+\Big(\sigma - \frac{p}{2} \Big ) \tau R_\alpha A \Big ) \frac{y_\alpha^{n+1} - y^n}{p \tau} +
R_\alpha A \frac{y_\alpha^{n+1} + y^n}{2} = 0
\]
for $\alpha=1,2,\ldots,p$.
Multiplying by $A^{1/2}$, we arrive at Eq.
\[
\frac{\widetilde{y}_\alpha^{n+1} - \widetilde{y}^n}{p \tau} +
\Big (I+\Big(\sigma - \frac{p}{2} \Big ) \tau \widetilde{A}_\alpha \Big )^{-1} \widetilde{A}_\alpha \frac{\widetilde{y}_\alpha^{n+1} + \widetilde{y}^n}{2} = 0 ,
\]
in which $\widetilde{y}_\alpha^{n+1} = A^{1/2}y_\alpha^{n+1}$, $\widetilde{y}^n = A^{1/2}y^n$.
Multiplying scalarly in $H$ by $p\tau(\widetilde{y}_\alpha^{n+1} + \widetilde{y}^n)$ and considering the non-negativity and symmetry of the operators
$ \widetilde{A}_\alpha$, we obtain
\[
\|\widetilde{y}_\alpha^{n+1}\| \leq \|\widetilde{y}^n \|,
\quad \alpha=1,2,\ldots,p .
\]
Given the representation (\ref{4.11}) for the solution at the new time level, we arrive at the provable estimate (\ref{4.9}).
\end{pf}

\subsection{Second decomposition variant}

Using the decomposition variant (\ref{4.3}) we have the equation
\begin{equation}\label{4.12}
\frac{d u}{d t} + \sum_{\alpha = 1}^{p} A R_\alpha u = f(t),
\quad 0 < t \leq T .
\end{equation}
To symmetrize it, we multiply the equation by $A^{-1/2}$ and obtain equation (\ref{4.5}),
in which now $\widetilde{u} = A^{-1/2} u, \ \widetilde{f} = A^{-1/2} f$.
The initial condition (\ref{2.2}) in the new notations is written in the form of
\begin{equation}\label{4.13}
\widetilde{u} (0) = \widetilde{u}^0,
\quad \widetilde{u}^0 = A^{-1/2} u^0 .
\end{equation}
For the solution, there is an estimate (\ref{2.4}) for $D = A^{-1}$.

We can carry out composition schemes for approximate solutions of the problem (\ref{2.2}), (\ref{4.1}) based on the transition to the problem (\ref{4.5}), (\ref{4.13}) with symmetrized operators $\widetilde{A}_\alpha, \ \alpha=1,2,\ldots,p$.
We discussed a similar approach when considering the first decomposition (\ref{4.2}).
Let us note another possibility of constructing composition schemes, which is related to considering an auxiliary problem for a system of evolution equations.
Here, we follow the work \cite{vabishchevich2023difference}, in which two- and three-level domain decomposition schemes are constructed.

Multiplying equation (\ref{4.12}) by $R_\alpha, \ \alpha = 1,2, \ldots, p$, we get
\begin{equation}\label{4.14}
R_\alpha \frac{d u}{d t} + R_\alpha \sum_{\beta =1}^{p} A R_\beta u = f_\alpha (t),
\quad 0 < t\leq T ,
\quad \alpha = 1,2, \ldots, p ,
\end{equation}
where now
\[
f_\alpha (t) = R_\alpha f(t),
\quad \alpha = 1,2, \ldots, p .
\]
Let's define auxiliary functions $u_\alpha, \ \alpha = 1,2, \ldots, p$ to be determined from Eqs.
\begin{equation}\label{4.15}
R_\alpha \frac{d u_\alpha }{d t} + R_\alpha \sum_{\beta =1}^{p} A R_\beta u_\beta = R_\alpha f(t),
\quad 0 < t\leq T ,
\quad \alpha = 1,2, \ldots, p .
\end{equation}
Given (\ref{4.15}), each function $u_\alpha, \ \alpha = 1,2, \ldots, p$ is mapped to $u$.
Let us supplement the system of equations (\ref{4.15}) with initial conditions
\begin{equation}\label{4.16}
u_\alpha(0) = u^0,
\quad \alpha = 1,2, \ldots, p .
\end{equation}

Let's define
\begin{equation}\label{4.17}
u(t) = \sum_{\alpha =1}^{p} R_\alpha u_\alpha(t) ,
\quad 0 < t\leq T .
\end{equation}
Adding the equations (\ref{4.15}), we obtain
\[
\sum_{\alpha =1}^{p} R_\alpha \frac{d u_\alpha }{d t} + \sum_{\alpha =1}^{p} R_\alpha A \sum_{\beta =1}^{p} R_\beta u_\beta =
\sum_{\alpha =1}^{p} R_\alpha f(t) .
\]
It follows that if $u_\alpha(t), \ \alpha = 1,2, \ldots, p$ is a solution to the Cauchy problem (\ref{4.15}), (\ref{4.16}),
then $u(t)$, defined according to (\ref{4.16}), is a solution of the problem (\ref{2.1}), (\ref{2.2}).
Thus, we can construct an approximate solution at a new time level in the form (\ref{4.17}) based on explicit-implicit approximations for the system of equations (\ref{4.15}).

Suppose, for example, that for (\ref{4.15}), (\ref{4.16}) with $f_\alpha(t) = 0, \ \alpha = 1,2, \ldots, p$ a two-level scheme is used
\begin{equation}\label{4.18}
R_\alpha \frac{y^{n+1}_\alpha - y^{n}_\alpha}{\tau} +
\sigma R_\alpha A R_\alpha (y^{n+1}_\alpha - y^{n}_\alpha) +
R_\alpha \sum_{\beta =1}^{p} A R_\beta y^{n}_\beta = 0,
\quad n = 0,1, \ldots, N-1 ,
\end{equation}
\begin{equation}\label{4.19}
y^{0}_\alpha = u^0 ,
\quad \alpha = 1,2, \ldots, p .
\end{equation}

\begin{thm}\label{t-3}
When $\sigma \geq p/2$ for the scheme (\ref{4.18}), (\ref{4.19}) there is an a priori estimate (\ref{4.9}) given by
\begin{equation}\label{4.20}
y^{n+1} = \sum_{\alpha = 1}^{p} R_\alpha y_\alpha^{n+1} ,
\quad n = 0,1, \ldots, N-1 .
\end{equation}
\end{thm}

\begin{pf}
Let us write (\ref{4.18}) in the form
\begin{equation}\label{4.21}
\Big ( R_\alpha + \sigma \tau R_\alpha A R_\alpha \Big ) \frac{y^{n+1}_\alpha - y^{n}_\alpha}{\tau} - \frac{\tau }{2} R_\alpha \sum_{\beta =1}^{p} A R_\beta \frac{y^{n+1}_\beta - y^{n}_\beta}{\tau}
+
R_\alpha \sum_{\beta =1}^{p} A R_\beta \frac{y^{n+1}_\beta + y^{n}_\beta}{2} = 0 .
\end{equation}
Let's multiply the individual equations scalarly in $H$ by
\[
v_\alpha = 2 (y^{n+1}_\alpha - y^{n}_\alpha)
\]
and add them up.
Taking into account (\ref{4.20}), for the last summand in (\ref{4.21}) we obtain
\[
\Big ( \sum_{\alpha =1}^{p} R_\alpha (y^{n+1}_\alpha - y^{n}_\alpha), A R_\alpha \sum_{\beta =1}^{p} (y^{n+1}_\beta + y^{n}_\beta) \Big ) = \big ((y^{n+1} - y^n), A (y^{n+1} + y^n) \big ) = \|y^{n+1}\|_A^2 - \|y^{n}\|_A^2 .
\]
Taking into account the introduced notations for other summands in (\ref{4.21}) we have
\[
\frac{2}{\tau } \sum_{\alpha =1}^{p} \big ( R_\alpha v_\alpha, v_\alpha) \geq 0 ,
\]
\[
2 \sigma \sum_{\alpha =1}^{p} ( R_\alpha v_\alpha, A R_\alpha v_\alpha) = 2 \sigma \sum_{\alpha =1}^{p} (w_\alpha,w_\alpha) ,
\]
\[
\sum_{\alpha =1}^{p} \Big ( R_\alpha v_\alpha, \sum_{\beta =1}^{p} A R_\beta v_\beta \beta \Big ) =
\sum_{\alpha =1}^{p} \Big (w_\alpha , \sum_{\beta =1}^{p} w_\beta \Big ) \leq p \sum_{\alpha =1}^{p} (w_\alpha,w_\alpha) ,
\]
where $w_\alpha = A^{1/2} R_\alpha v_\alpha, \ \alpha = 1,2, \ldots, p$.
With the noted constraints on $\sigma$, we obtain the inequality
\[
\|y^{n+1}\|_A^2 \leq \|y^{n}\|_A^2 ,
\]
which gives the provable estimate (\ref{4.9}).
\end{pf}

Explicit-implicit time approximations are similarly considered when instead of (\ref{4.18}), one uses
\begin{equation}\label{4.22}
R_\alpha \frac{y^{n+1}_\alpha - y^{n}}{\tau} +
\sigma R_\alpha A R_\alpha (y^{n+1}_\alpha - y^{n}) +
R_\alpha \sum_{\beta =1}^{p} A R_\beta y^{n} = 0,
\quad n = 0,1, \ldots, N-1 ,
\end{equation}
with the composition of the approximate solution according to (\ref{4.20}) at each level in time.
The resulting scheme (\ref{2.6}), (\ref{4.20}), (\ref{4.22}) can be related to special variants of the additive-averaging scheme.

\subsection{Other spaces}

We constructed composition schemes for the problem (\ref{2.1}), (\ref{2.2}) based on the additive representation (\ref{4.1}) of the unit operator in space $H$.
One can be guided to work in a different space by attracting additional information about the equation.
In particular, this approach allows us to consider classical splitting schemes on spatial variables for the approximate solution of the model problem (\ref{3.1})--(\ref{3.3}).

Consider the equation (\ref{2.1}) in which the operator $A \geq 0$ is factorized:
\begin{equation}\label{4.23}
A = D^* D ,
\end{equation}
and $D: H \rightarrow \bm H$. Now the decomposition for the unit operator $\bm I$ in $\bm H$ is used:
\begin{equation}\label{4.24}
\bm I = \sum_{\alpha = 1}^{p} R_\alpha ,
\quad R_\alpha = R_\alpha^* \geq 0,
\quad \alpha=1,2,\ldots,p .
\end{equation}

Given (\ref{4.22}), (\ref{4.23}), we get the equation
\begin{equation}\label{4.25}
\frac{d u}{d t} + \sum_{\alpha = 1}^{p} D^* R_\alpha D u = f(t),
\quad 0 < t \leq T .
\end{equation}
In the case of (\ref{4.25}), we have a decomposition of the problem operator
\[
A = \sum_{\alpha =1}^{p} A_\alpha \geq 0 ,
\quad A_\alpha = D^* R_\alpha D = A_\alpha^* \geq 0,
\quad \alpha = 1,2, \ldots, p .
\]
As such, there are no problems constructing composition schemes based on different classes of additive operator schemes \cite{VabishchevichAdditive}.
Stability (a priori estimates of the type (\ref{2.7})) for the approximate solution takes place in $H$.

Let us illustrate the decomposition (\ref{4.23}), (\ref{4.24}) on the example of the problem (\ref{3.1})--(\ref{3.3}), when the operator $D$ is associated with the gradient operator and the operator $D^*$ --- with the divergence operator.
Let us define the Hilbert space $\bm H$ of vector grid functions $\bm v = \{v_1(\bm x), v_2(\bm x)\}, \ \bm x \in \omega$, in which
\[
(\bm v, \bm w) = (v_1, w_1) + (v_2, w_2),
\quad \|\bm v\| = (\bm v, \bm v)^{1/2} .
\]
For the grid operator $A$, defined according to (\ref{3.7}), for $D$ we obtain
\[
\begin{split}
D u = \{D_1 u, D_2 u\},
\quad & D_1 u = \frac{1}{h_1} k^{1/2}(x_1+0.5h_1,x_2)
\big (u(x_1+h_1,x_2) - u(x_1,x_2) \big ) , \\
\quad & D_2 u = \frac{1}{h_2} k^{1/2}(x_1,x_2+0.5h_2)
\big (u(x_1,x_2+h_2) - u(x_1,x_2) \big ) .
\end{split}
\]
The splitting schemes on spatial variables correspond to the assignment
\[
R_1 \bm v = \begin{pmatrix}
v_1(\bm x) & 0 \\
0 & 0
\end{pmatrix} ,
\quad R_2 \bm v = \begin{pmatrix}
0 & 0\\
0 & v_2(\bm x)
\end{pmatrix} ,
\quad \bm x \in \omega .
\]
When constructing domain decomposition schemes, we take the additive representation (\ref{4.24}) to be
\[
R_\alpha \bm v = \begin{pmatrix}
\chi_\alpha (\bm x) v_1(\bm x) & 0 \\
0 &\chi_\alpha (\bm x) v_2(\bm x)
\end{pmatrix} ,
\quad \alpha = 1,2, \ldots, p ,
\quad \bm x \in \omega .
\]
Here, we used the partitioning of the set of internal grid nodes (\ref{3.15}).

\section{Decompositions for a set of spaces}\label{sec:5}

The computational technique of decomposition and composition can be realized by defining operators $R_\alpha, \ \alpha = 1,2, \ldots, p$ in various spaces.
In our study, we are based on the work of \cite{vabishchevich2023subdomain}, in which such approximate solution constructions are realized as domain decomposition schemes.

\subsection{Additive representation of the unit operator}

In addition to the space $H$, our consideration includes a set of finite-dimensional Hilbert spaces $H_\alpha, \ \alpha = 1,2,\ldots, p$.
For a particular space, we define the linear operators (the restriction operators) $G_\alpha$:
\[
G_\alpha: H \rightarrow H_\alpha,
\quad G^*\alpha: H\alpha \rightarrow H,
\quad \alpha = 1,2,\ldots,p .
\]
These operators define the decomposition of the unit operator in $H$ by the rule
\begin{equation}\label{5.1}
I = \sum_{\alpha =1}^{p} G_\alpha^* G_\alpha .
\end{equation}
This construction was previously used in the works of \cite{vabishchevich2021solution,vabishchevich2017vector}.
When decomposing the computational grid according to (\ref{3.15}) we have (see \cite{vabishchevich2023subdomain})
\[
G_\alpha u = \chi^{1/2}_\alpha (\bm x) u,
\quad \bm x \in \omega ,
\quad \alpha=1,2,\ldots,p .
\]

The decomposition (\ref{5.1}) corresponds to the use of (\ref{4.1}) at
\[
R_\alpha = G_\alpha^* G_\alpha ,
\quad \alpha = 1,2,\ldots,p .
\]
Given (\ref{5.1}), we can directly obtain an additive representation for the solution $u \in H$ via the auxiliary solutions $u_\alpha \in H_\alpha, \ \alpha = 1,2,\ldots,p$:
\begin{equation}\label{5.2}
u = \sum_{\alpha =1}^{p} G_\alpha^* G_\alpha u = \sum_{\alpha=1}^{p} G_\alpha^* u_\alpha,
\quad u_\alpha = G_\alpha u,
\quad \alpha = 1,2,\ldots,p .
\end{equation}
Thus, from the decomposition of the unit operator (\ref{5.1}), the composition of the solution (\ref{5.2}) immediately follows.

\subsection{Problem for solution components}

When solving (\ref{2.1})--(\ref{2.3}), the decomposition of (\ref{5.1}) gives Eq.
\[
\frac{d u}{d t} + \sum_{\beta = 1}^{p} A G_\beta^* G_\beta \, u = f(t),
\quad 0 < t \leq T .
\]
Given (\ref{5.2}), we obtain
\begin{equation}\label{5.3}
\frac{d u}{d t} + \sum_{\beta = 1}^{p} A G_\beta^* \, u_\beta = f(t),
\quad 0 < t \leq T .
\end{equation}
Let's multiply (\ref{5.3}) by $G_\alpha, \ \alpha = 1,2,\ldots, p$, which gives a system of equations
\begin{equation}\label{5.4}
\frac{d u_\alpha }{d t} + \sum_{\beta = 1}^{p} G_\alpha A G_\beta^* \, u_\beta = G_\alpha f(t),
\quad \alpha=1,2,\ldots,p ,
\quad 0 < t \leq T ,
\end{equation}
for the individual components of the solution $u_\alpha \in H_\alpha , \ \alpha = 1,2,\ldots, p$.
Considering (\ref{2.2}), we augment this system of evolution equations with initial conditions
\begin{equation}\label{5.5}
u_\alpha(0) = u_\alpha^0,
\quad u_\alpha^0 = G_\alpha u^0 ,
\quad \alpha = 1,2,\ldots,p .
\end{equation}
For $u \in H$, we have a representation
\begin{equation}\label{5.6}
u = \sum_{\alpha=1}^{p} G_\alpha^* u_\alpha .
\end{equation}

\begin{thm}\label{t-4}
Let $u_\alpha(t) \in H_\alpha , \ \alpha = 1,2,\ldots, p$ be the solution of the problem (\ref{5.4}), (\ref{5.5}).
Then $u(t)$ defined according to (\ref{5.6}) is a solution to the Cauchy problem (\ref{2.1}), (\ref{2.2}).
\end{thm}

\begin{pf}
Under condition (\ref{5.6}), equations (\ref{5.4}) are written as
\[
\frac{d u_\alpha }{d t} + G_\alpha A u = G_\alpha f(t),
\quad \alpha=1,2,\ldots,p .
\]
Multiplying the individual equation by $G_\alpha^*$ and adding them together, given (\ref{5.1}) and (\ref{5.6}), we get the equation (\ref{2.1}).
Similarly, from (\ref{5.5}) follows (\ref{2.2}).
\end{pf}

\subsection{Time approximation}

Various explicit-implicit schemes can be used to approximate the solution of the problem (\ref{5.4})--(\ref{5.6}).
Here is a typical result.
For the Cauchy problem (\ref{5.4}), (\ref{5.5}) with $f(t) = 0$ we will use a two-level scheme
\begin{equation}\label{5.7}
\frac{y^{n+1}_\alpha - y^{n}_\alpha}{\tau} +
\sigma G_\alpha A G^*_\alpha (y^{n+1}_\alpha - y^{n}_\alpha) +
G_\alpha \sum_{\beta =1}^{p} A G^*_\beta y^{n}_\beta = 0,
\quad n = 0,1, \ldots, N-1 ,
\end{equation}
\begin{equation}\label{5.8}
y^{0}_\alpha = u_\alpha^0 ,
\quad \alpha = 1,2, \ldots, p .
\end{equation}
Similarly to Theorem~\ref{t-3}, the following statement is proved.

\begin{thm}\label{t-6}
When $\sigma \geq p/2$ for the scheme (\ref{5.7}), (\ref{5.8}) there is an a priori estimate (\ref{4.9}) on the composition of the solution
\begin{equation}\label{5.9}
y^{n+1} = \sum_{\alpha = 1}^{p} G^*_\alpha y_\alpha^{n+1} ,
\quad n = 0,1, \ldots, N-1 .
\end{equation}
\end{thm}

\begin{pf}
Let's rewrite (\ref{5.7}) in the form of
\begin{equation}\label{5.10}
\Big ( I + \sigma \tau G_\alpha A G^*_\alpha \Big ) \frac{y^{n+1}_\alpha - y^{n}_\alpha}{\tau} - \frac{\tau }{2} G_\alpha \sum_{\beta =1}^{p} A G^*_\beta \frac{y^{n+1}_\beta - y^{n}_\beta}{\tau}
+ G_\alpha \sum_{\beta =1}^{p} A G^*_\beta \frac{y^{n+1}_\beta + y^{n}_\beta}{2} = 0 .
\end{equation}
Let's multiply the individual equations scalarly in $H$ by
\[
v_\alpha = 2 (y^{n+1}_\alpha - y^{n}_\alpha)
\]
and add them up.
For the individual summands in (\ref{5.10}), this gives
\[
2 \sigma \sum_{\alpha =1}^{p} ( G^*_\alpha v_\alpha, A G^*_\alpha v_\alpha) = 2 \sigma \sum_{\alpha =1}^{p} (w_\alpha,w_\alpha) ,
\]
\[
\sum_{\alpha =1}^{p} \Big ( G^*_\alpha v_\alpha, \sum_{\beta =1}^{p} A G^*_\beta v_\beta \Big ) =
\sum_{\alpha =1}^{p} \Big (w_\alpha , \sum_{\beta =1}^{p} w_\beta \Big ) \leq p \sum_{\alpha =1}^{p} (w_\alpha,w_\alpha) ,
\]
\[
\Big ( \sum_{\alpha =1}^{p} G^*_\alpha (y^{n+1}_\alpha - y^{n}_\alpha), A G^*_\alpha \sum_{\beta =1}^{p} (y^{n+1}_\beta + y^{n}_\beta) \Big ) = \big ((y^{n+1} - y^n), A (y^{n+1} + y^n) \big ) = \|y^{n+1}\|_A^2 - \|y^{n}\|_A^2 .
\]
where $w_\alpha = A^{1/2} G^*_\alpha v_\alpha, \ \alpha = 1,2, \ldots, p$.
When $\sigma \geq p/2$, we obtain the estimate (\ref{4.9}).
\end{pf}

Some other possibilities for constructing two- and three-level composition schemes can be found in \cite{vabishchevich2023subdomain}.
In particular, let us note the three-level scheme
\[
\frac{y^{n+1}_\alpha - y^{n-1}_\alpha}{2\tau} +
\sigma G_\alpha A G^*_\alpha (y^{n+1}_\alpha - 2y^{n}_\alpha + y^{n-1}_\alpha) +
G_\alpha \sum_{\beta =1}^{p} A G^*_\beta y^{n}_\beta = 0,
\quad n = 0,1, \ldots, N-1 ,
\]
with determination of the approximate solution according to (\ref{5.9}), which has a second-order approximation on $\tau$ and is stable at $\sigma \geq p/4$.

\section{Other problems}\label{sec:6}

We have considered the Cauchy problem (\ref{2.1})--(\ref{2.3}) for a first-order evolutionary with a self-adjoint positive operator.
We mention the possibilities of using the proposed decomposition and composition technique to solve other nonstationary problems.

\subsection{Problems with non-self-adjoint operators}

Let us assume that in Eq. (\ref{2.1})
\[
A = B + C,
\quad B = \frac{1}{2} (A + A^* ) \geq 0,
\quad C = - C^* = \frac{1}{2} (A - A^* ) ,
\]
so that $A \geq 0$.
For the self-adjoint part of the operator $A$, we can use the decomposition
\[
B = \sum_{\alpha = 1}^{p} B_\alpha,
\quad B_\alpha = B_\alpha^* \geq 0,
\quad \alpha = 1,2, \ldots, p ,
\]
based on the variant decomposition of the unit operator (\ref{4.24}).
The problem of decomposition of the skew-symmetric part of the operator $A$ deserves special attention.

When using the additive representation of the unit operator (\ref{4.1}), the following is true
\[
C = \frac{1}{2} \sum_{\alpha = 1}^{p} \big ( R_\alpha C + C R_\alpha \big ) .
\]
By this we have the decomposition
\[
C = \sum_{\alpha = 1}^{p} C_\alpha,
\quad C_\alpha = - C_\alpha^* = \frac{1}{2} \big ( R_\alpha C + C R_\alpha \big ) ,
\quad \alpha = 1,2, \ldots, p .
\]
We have constructed (see, e.g., \cite{Vabischevich1996c,SamarskiiVabischevich1997a}) domain decomposition schemes for unsteady convection-diffusion problems along such a path.
For some operator summands, we use one class of decomposition (decomposition of the unit operator in one space), and for other operator summands --- another class of decomposition (decomposition of the unit operator in another space).

The composition schemes are more easily constructed when working with a set of spaces.
When the unit operator is decomposed in the form (\ref{5.1}), the basic system of equations is (\ref{5.4}), in which the critical property of non-negativity of the non-self-adjoint operator $A$ is inherited by the individual operators $A_\alpha, \ \alpha = 1,2, \ldots, p$:
\[
A = \sum_{\alpha =1}^{p} A_\alpha \geq 0 ,
\quad A_\alpha = G_\alpha A G_\alpha^* \geq 0,
\quad \alpha = 1,2, \ldots, p .
\]
When investigating the time approximation of the corresponding composition schemes, we should focus on obtaining stability estimates in $H$.

\subsection{Evolutionary equations of second order}

It is fundamentally possible to construct decomposition-composition schemes for second-order evolutionary equations.
We will find an approximate solution to the Cauchy problem
\begin{equation}\label{6.1}
\frac{d^2 u}{d t^2} + A u = f(t),
\quad 0 < t \leq T,
\end{equation}
\begin{equation}\label{6.2}
u(0)= u^0 ,
\quad \frac{d u}{d t} (0) = v^0 ,
\end{equation}
under the assumptions (\ref{2.3}).

The decomposition of the problem (\ref{6.1}), (\ref{6.2}) follows the same scheme as for the first-order evolution equations.
For example, for the additive representation of the unit operator according to (\ref{4.1}), the decomposition in (\ref{4.2}) gives Eq.
\begin{equation}\label{6.3}
\frac{d^2 u}{d t^2} + \sum_{\alpha = 1}^{p} R_\alpha A u = f(t),
\quad 0 < t \leq T ,
\end{equation}
which is symmetrized in the usual way.

In composing the solution, we apply the results of the theory of additive operator-difference schemes \cite{VabishchevichAdditive} and construct various explicit-implicit three-level time approximations.
For example, similar to (\ref{4.8}), a regularized scheme can be used for Eq. (\ref{6.3})
\[
\frac{y^{n+1} - 2 y^n + y^{n-1}}{\tau^2} + \sum_{\alpha = 1}^{p} (I+\sigma \tau^2 R_\alpha A)^{-1} R_\alpha A y^n = f^{n} .
\]
It is unconditionally stable at $\sigma \geq p/4$.

\subsection{Systems of equations}

An important promising direction in the theory and practice of the computational technology of decomposition and composition is its application in the approximate solution of nonstationary problems for systems of evolutionary equations.
For such problems, we construct computational algorithms when the transition to a new level in time is carried out by solving more straightforward issues for the individual quantities required.
As a model problem, we consider the Cauchy problem for a system of two equations.

We find $u_\alpha(t) \in H, \ \alpha =1,2$ from the equations.
\begin{equation}\label{6.4}
\begin{split}
\frac{d u_1}{d t} + A_{11} u_1 + A_{12} u_2 & = f_1(t), \\
\frac{d u_1}{d t} + A_{21} u_1 + A_{22} u_2 & = f_2(t),
\quad 0 < t \leq T ,
\end{split}
\end{equation}
with constant operators $A_{\alpha \beta}: H \rightarrow H, \ \alpha, \beta = 1,2$
and the initial conditions
\begin{equation}\label{6.5}
u_1(0) = u_1^0,
\quad u_2(0) = u_2^0 .
\end{equation}
Difficulties in the numerical solution of this problem may be because the system of equations is coupled due to operators $A_{12}, \ A_{21}$.

Let us introduce the vector of sought quantities $\bm u = \{u_1, u_2\}$, for the initial data and the right-hand side we put $\bm u^0 = \{u_1^0, u_2^0\}$, $\bm f = \{f_1, f_2\}$, respectively.
Let us write the system of equations (\ref{6.4}) in the form of
\begin{equation}\label{6.6}
\frac{d \bm u}{d t} + \bm A \bm u = \bm f ,
\quad 0 < t\leq T .
\end{equation}
For the operator matrix $\bm A$ we have
\begin{equation}\label{6.7}
\quad \bm A = \begin{pmatrix}
A_{11} & A_{12} \\
A_{21} & A_{22}
\end{pmatrix} .
\end{equation}
Equation (\ref{6.6}) is supplemented by the initial condition
\begin{equation}\label{6.8}
\bm u(0) = \bm u^0 .
\end{equation}

We will consider the problem (\ref{6.6})--(\ref{6.8}) on the direct sum of spaces $\bm H = H \oplus H$.
For $\bm u, \bm v \in \bm H$, the scalar product and norm are
\[
(\bm u, \bm v) = \sum_{\alpha =1}^{2} (u_\alpha, v_\alpha) ,
\quad \|\bm u\| = (\bm u, \bm u)^{1/2} .
\]
Let $A_{\alpha \alpha} = A_{\alpha \alpha}^*, \ \alpha = 1,2$ and $A_{21} = A_{12}^*$ in $H$ and so $\bm A = \bm A^*$ in $\bm H$.
Similarly (\ref{2.3}), we will assume that the operator $\bm A$ is positive.
If $A_{22}^{-1}$ exists, this property requires the operator inequalities
\[
A_{22} > 0,
\quad A_{11} > A_{12} A_{22}^{-1} A_{12}^* .
\]

We will use the decomposition of the unit operator $\bm I$ into $\bm H$ in the form (see (\ref{4.1}))
\[
\bm I = \bm R_1 + \bm R_2,
\quad \bm R_\alpha = \bm R_\alpha^* \geq 0,
\quad \alpha = 1,2 ,
\]
in which
\[
\bm R_1 = \begin{pmatrix}
I & 0 \\
0 & 0
\end{pmatrix} ,
\quad \bm R_2 = \begin{pmatrix}
0 & 0 \\
0 & I
\end{pmatrix} .
\]
In the decomposition variant (\ref{4.2}), we get the equation
\begin{equation}\label{6.9}
\frac{d \bm u}{d t} + \bm A_1 \bm u + \bm A_2 \bm u = \bm f ,
\quad 0 < t\leq T ,
\end{equation}
in which
\[
\bm A_1 = \begin{pmatrix}
A_{11} & A_{12} \\
0 & 0
\end{pmatrix} ,
\quad \bm A_2 = \begin{pmatrix}
0 & 0 \\
A_{21} & A_{22}
\end{pmatrix} .
\]
A variant decomposition (\ref{4.2}) yields
\begin{equation}\label{6.10}
\bm A_1 = \begin{pmatrix}
A_{11} & 0 \\
A_{21} & 0
\end{pmatrix} ,
\quad \bm A_2 = \begin{pmatrix}
0 & A_{12} \\
0 & A_{22}
\end{pmatrix} .
\end{equation}
We have a split of the operator matrix $\bm A$ either by rows or by columns.

When applying explicit-implicit approximations for equation (\ref{6.9}), we obtain separate problems for the solution components $u_\alpha(t) \in H, \ \alpha =1,2$ at a new level in time.
For example, using a purely implicit component-wise splitting scheme (see (\ref{3.22}) at $\sigma = 1$) gives at $\bm f(t) = 0$
\[
\frac{\bm y^{n+1/2} - \bm y^{n}}{\tau } + \bm A_1 \bm y^{n+1/2} = 0,
\]
\[
\frac{\bm y^{n+1} - \bm y^{n+1/2}}{\tau } + \bm A_2 \bm y^{n+1} = 0 .
\]
In the case of decomposition (\ref{6.10}), we obtain
\[
\frac{y_1^{n+1/2} - y_1^{n}}{\tau } + A_{11} y_1^{n+1/2} = 0,
\quad \frac{y_2^{n+1/2} - y_2^{n+1/2}}{\tau } + A_{21} y_1^{n+1/2} = 0,
\]
\[
\frac{y_1^{n+1} - y_1^{n+1/2}}{\tau } + A_{12} y_2^{n+1} = 0 ,
\quad \frac{y_2^{n+1} - y_2^{n+1/2}}{\tau } + A_{22} y_2^{n+1} = 0 .
\]
The calculations are performed in the following sequence: $y_1^{n+1/2}$, $y_2^{n+1/2}$, $y_2^{n+1}$, $y_1^{n+1}$.
We solve problems only with diagonal elements $A_{11}, A_{22}$ of the operator matrix $\bm A$.

\section{Conclusions} {\label{sec:7}

\begin{enumerate}[(1)]

\item Based on the practice and theory of construction and research of additive operator schemes (splitting schemes), this paper formulates a computational technique of decomposition and composition for approximate solution of Cauchy problems for evolution equations, which are considered in a finite-dimensional Hilbert space.
The problem operator is a sum of more straightforward operators in the decomposition stage.
At the composition stage, an approximate solution to the problem is constructed based on the solution of auxiliary problems using various explicit-implicit time approximations.

\item The primary attention is paid to separating more straightforward auxiliary problems.
In the first variant, decomposition is provided by using the decomposition of the unit operator as a sum of self-adjoint non-negative operators. In the second decomposition variant, the individual operator summands in the unit operator decomposition are associated with different spaces.

\item Composition of the approximate solution is provided by applying different variants of additive operator-difference schemes.
For the auxiliary problems we obtain at the decomposition stage, splitting schemes of different classes are constructed.
The study focuses on obtaining stability estimates of the approximate solution in the corresponding space.

\item The main elements of the decomposition-composition technique are illustrated on the Cauchy problem for the first-order evolution equation with a self-adjoint positive operator.
After applying the standard finite-difference approximation over space for the second-order parabolic equation in a rectangle, the differential-difference system of equations is a test problem. The decomposition is related to the splitting over spatial variables and the decomposition of the computational grid.

\item The possibilities of application of the decomposition-composition technique to other problems are briefly mentioned: problems with non-self-adjoint operators and problems for second-order evolutionary equations.
We currently associate the most interesting new applications of the developed general approach with systems of equations.
The main problems are illustrated in the Cauchy problem for a system of two coupled evolution equations of first order.

\end{enumerate}

\end{document}